\documentclass[11pt]{amsart}
\usepackage{latexsym,amssymb,amsmath}
\newtheorem{theorem}{Theorem}[section]
\newtheorem{proposition}[theorem]{Proposition}

\newtheorem{conjecture}[theorem]{Conjecture}

\theoremstyle{definition}

\newtheorem{exercise}[theorem]{Exercise}

\newtheorem{example}[theorem]{Example}

\theoremstyle{remark}
\newtheorem{remark}[theorem]{Remark}

\numberwithin{equation}{section}

\def\t{\tau}
\def\s{\sigma}

\def\frp#1{\frac{\partial}{\partial{#1}}}

\def\na{n+\alpha}

\def\xx#1{x^{#1}}

\def\ooo#1#2{\omega^{#1}_{#2}}\def\xxx#1#2{x^{#1}_{#2}}
\def\oo#1{\omega^{#1}_0}
\def\ep{\epsilon}
\def\qq#1#2#3{q^{#1}_{{#2} {#3}}}
\def\rr#1#2#3#4{r^{#1}_{{#2} {#3}{#4}}}
\def\rrr#1#2#3#4#5{r^{#1}_{{#2} {#3}{#4}{#5}}}
\def\ee#1{e_{#1}}

\def\ue#1{e^{#1}}
\def\cf{\mathcal F}

\def\trank{\text{rank}}

\def\BC{\mathbb C}\def\BF{\mathbb F}\def\BO{\mathbb O}
\def\BR{\mathbb R}
\def\BG{\mathbb G}
\def\BP{\mathbb P}
\def\pp#1{\mathbb P^{#1}}
\def\fa{\mathfrak a}\def\fc{\mathfrak c}
\def\fgl{\mathfrak g\mathfrak l}

\def\pp#1{{\mathbb P}^{#1}}
\def\tdim{\rm dim}
\def\hd{,...,}
\def\ww{\wedge}
\def\upperp{{}^\perp}

\def\na{n+a}

\def\inv{{}^{-1}}

\def\cC{{\mathcal C}}

\def\BZ{\mathbb Z}

\def\11{\mathbf 1}

\def\FF{\mathbb F}

\def\fsl{{\mathfrak {sl}}}

\def\fe{{\mathfrak e}}

\def\fg{{\mathfrak g}}

\def\fp{{\mathfrak p}}

\def\l{\lambda}
\def\a{\alpha}

\def\b{\beta}
\def\g{\gamma}
\def\s{\sigma}

\def\up#1{{}^{({#1})}}

\def\ot{{\mathord{\,\otimes }\,}}
\def\op{{\mathord{\,\oplus }\,}}

\def\ra{{\mathord{\;\rightarrow\;}}}

\def\La#1{\Lambda^{#1}}
\def\tim{\text{Image}\,}
\def\tdeg{\text{deg} }

\def\tdim{\text{dim}\,}
\def\tcodim{\text{codim}\,}

\def\tker{\text{ker}\,}

\def\tmod{\text{ mod }}
\def\tmin{\text{ min }}

\def\trank{\text{rank}\,}

\begin{document}
\title{Differential geometry of submanifolds of projective space}
\author{J.M. Landsberg }
\date{September 2006}
\begin{abstract} These are lecture notes on the   rigidity
of submanifolds of projective space \lq\lq resembling\rq\rq\ 
compact Hermitian symmetric spaces in their homogeneous embeddings.
The results of \cite{Lrigid,Lchss,LR,Lline,Lnfac,HY,robles} are
surveyed, along with their classical predecessors. The notes
include an introduction to moving frames in projective geometry,
an exposition of the Hwang-Yamaguchi ridgidity theorem and
  a new variant of the Hwang-Yamaguchi theorem.
\end{abstract}
\thanks{Supported by NSF grant DMS-0305829}
\email{jml@math.tamu.edu}
\maketitle

    

\section{Overview}
\begin{itemize}
\item Introduction to the local differential geometry of submanifolds of projective space.
\item Introduction to   moving frames for projective geometry.
\item   How much must a submanifold $X\subset \BP^N$ resemble 
a given submanifold $Z\subset \BP^M$ infinitesimally before
we can conclude $X\simeq Z$?
\item   To what order must a line field on a submanifold $X\subset \BP^N$ have
contact with $X$ before we can conclude the lines are contained in $X$?
\item Applications to algebraic geometry.
\item A new variant of the Hwang-Yamaguchi rigidity theorem.
\item An exposition of the Hwang-Yamaguchi rigidity theorem in the language
of moving frames.
\end{itemize}

Representation theory and algebraic geometry
are natural tools for studying submanifolds of projective
space. Recently  there has also been progress
the other way, using projective differential geometry
to prove results in algebraic geometry and representation theory. These talks will focus
on the basics of submanifolds of projective
space, and give a few applications to
algebraic geometry. For further applications
to algebraic geometry the reader is invited
to consult chapter 3 of \cite{EDSBK} and the references therein.

Due to constraints of time and space, applications
to representation theory will not be given here,
but the interested reader can consult
\cite{LMpop} for an overview. Entertaining applications
include new proofs of the classification of compact Hermitian
symmetric spaces, and of complex simple Lie algebras, based
on the geometry of rational homogeneous varieties (instead of root systems), see \cite{LM1}.
The  applications are not limited to classical
representation theory.  
There are applications to Deligne's conjectured categorical
generalization of the exceptional series \cite{LMdel}, to Vogel's proposed
{\it Universal Lie algebra} \cite{LMuniv}, and 
to the study of the intermediate Lie algebra $\fe_{7\frac 12}$ \cite{LMsex}.

\subsection*{Notations, conventions} I  mostly work over the
complex numbers in the complex analytic category, although
most of the results are valid in the $C^{\infty}$ category and over other fields, even
characteristic $p$, as long as the usual precautions are
taken. When working over $\BR$, some results become more complicated
as there are more possible normal forms. 
  I use notations and
the ordering of roots as in \cite{bou} and label maximal parabolic
subgroups accordingly, e.g., $P_k$ refers to the maximal parabolic
obtained by omitting the spaces corresponding to the simple
root $\a_k$. $\langle v_1\hd v_k\rangle$ denotes the linear span
of the vectors $v_1\hd v_k$. If $X\subset \BP V$ is a subset,
$\hat X\subset V$ denotes the corresponding cone in $V\backslash 0$,
the inverse image of $X$ under the projection $V\backslash 0\ra \BP V$.
and $\overline X\subset \BP V$ denotes the Zariski closure of $X$,
the zero set of all the homogeneous polynomials vanishing on $X$.
When we write $X^n$, we mean $\tdim (X)=n$.
We often use $Id$ to denote the identity matrix or identity map.
Repeated indicies are to be summed over.

\subsection*{Acknowledgements} These notes are based on lectures
given at Seoul National University  in June 2006, the IMA
workshop {\it Symmetries and overdetermined systems of partial
differential equations}, July 2006 and at CIMAT (Guanajuato) August 2006. It is a pleasure to thank Professors
Han, Eastwood and Hernandez for inviting me to give these respective lecture
series. I would also like to thank Professor Yamaguchi for carefully
explaining his results with Hwang to me at the workshop and Professor
  Robles for   reading a draft of this article and providing
  corrections and suggestions for improvement.

\section{Submanifolds of projective space}

\subsection{Projective geometry}
Let $V$ be a vector space and let  $\BP V$ denote
the associated projective space. We   think of $\BP V$ as the
quotient of $GL(V)$, the general linear group of
invertible endomorphisms of $V$, by the subgroup $P_1$
preserving a line.
For example if we take the line
$$
\langle \begin{pmatrix} 1\\ 0\\ \vdots\\ 0\end{pmatrix}\rangle
$$
then
$$
P_1=\begin{pmatrix} * & * &\hdots &*\\
0 &*&\hdots &*\\
\vdots &*&\hdots &*\\
0 &*&\hdots &*\end{pmatrix}
$$
where, if $\tdim V=N+1$ the blocking is $(1,N)\times (1,N)$,
so $P_1$ is the group of   invertible matrices with zeros
in the lower left hand block.

In the spirit of Klein, we consider two submanifolds
$M_1,M_2\subset \BP V$ to be  {\it equivalent} if there exists
some $g\in GL(V)$ such that $g.M_1=M_2$ and define the corresponding
notion of local equivalence.

Just as in the geometry of submanifolds of Euclidean
space, we will look for   {\it differential invariants}
that will enable us to determine if a neighborhood (germ) of
a point of $M_1$ is equivalent to a neighborhood (germ) of
a point of $M_2$. These invariants will be obtained by
taking derivatives at a point in a geometrically meaningful way.
Recall that   second derivatives  furnish  a complete set of differential
invariants 
 for surfaces in Euclidean
three space- the vector-bundle valued Euclidean first and
second fundamental forms, and two surfaces are locally
equivalent iff there exists a local diffeomorphism
$f: M_1\ra M_2$  preserving
the first and second fundamental forms.

The group of admissible motions  in projective space is  larger than   the
corresponding
Euclidean group so we expect to have to take more derivatives
to determine equivalence in the projective case than the Euclidean. For example, it
was long known that for hypersurfaces $X^n\subset\pp{n+1}$
that one needs at least three derivatives,
and Jensen and Musso \cite{JM}
showed that for most hypersurfaces, when $n>1$, three
derivatives are sufficient. For curves in the plane,
by a classical result of Monge, one
needs six derivatives!

In order to take derivatives in a way that will facilitate
extracting geometric information from them, we will use
the {\it moving frame}. Before developing the moving frame
in \S\ref{movingframesect}, we discuss
a few coarse invariants without machinery and state several
rigidity results.

\subsection{Asymptotic directions}
Fix $x\in X^n\subset \BP V$. After taking one derivative,
we have the tangent space $T_xX\subset T_x\BP V$, which is
the set of tangent directions to lines in $\BP V$ having
contact with $X$ at $x$ to order at least one. Since we
are discussing directions, it is better to consider
$\BP T_xX\subset \BP T_x\BP V$. Inside $\BP T_xX$ is
$\cC_{2,X,x}\subset \BP T_xX$, the set of tangent directions
to lines having contact at least two with $X$ at $x$, these
are called the {\it asymptotic directions } in Euclidean geometry,
and we continue to use the same terminology in the projective setting.
Continuing, we define $\cC_{k,X,x}$ for all $k$, and finally,
$\cC_{\infty,X,x}$, which, in the analytic category, equals $\cC_{X,x}$, the lines on (the completion of)
$X$ through $x$.
When $X$ is understood we sometimes write $\cC_{k,x}$ for $\cC_{k,X,x}$.
 

\medskip

{\it What does $\cC_{2,X,x}$, or more generally
$\cC_{k,X,x}$ tell us about the
geometry of $X$?}

\medskip

That is, what can we learn of the macroscopic geometry
of $X$ from the microscopic geometry at a point?
To increase the chances of getting meaningful information,
from now on, when we are in the analytic or algebraic
category,  we will work at a {\it general point}. Loosely
speaking, after taking $k$ derivatives there will be both
discrete   and continuous invariants. A general point is one where all the discrete invariants
are locally constant.

(To be more precise, if one is in the analytic category, one
should really speak of $k$-general points (those that
are general to order $k$), to insure there is
just a finite number of discrete invariants. In everything
that follows we will be taking just a finite number of derivatives
and we should say we are working at a $k$-general point where
$k$ is larger than the number of derivatives we are taking.)

When we are in the $C^{\infty}$ category, we will work in open
subsets and require whatever property we are studying at a point
holds at all points in the open subset.

For example, if $X^n$ is a hypersurface, then  $\cC_{2,X,x}$ is
a degree two hypersurface in $\BP T_xX$ (we will prove this below),
and thus its only invariant is its rank $r$. In particular, if $X$
is a smooth algebraic variety and $x\in X_{general}$ the rank is $n$ (see
e.g.,  \cite{GH,
EDSBK}) and thus we do not get much information. (In contrast, if
  $r<n$, then the Gauss map of $X$ is degenerate and $X$ is (locally)
ruled by $\pp{n-r}$'s.)

More generally, if $X^n\subset \pp{n+a}$, then
$\cC_{2,X,x}$ is
the intersection of at most $\tmin(a,\binom{n+1}2)$ quadric hypersurfaces, and one
generally expects that equality holds. In particular, if the codimension is
sufficiently large we expect $\cC_{2,x}$ to be empty and otherwise
it should have codimension $a$. When this fails to happen, there
are often interesting consequences for the macroscopic geometry
of $X$.

\subsection{The Segre variety and Griffiths-Harris conjecture}\label{segex}
Let $A,B$ be vector spaces and let $V=A\ot B$. Let
$$
X=\BP ({\rm rank\ one\ tensors})\subset \BP V.
$$
Recall that every rank one matrix
(i.e., rank one tensor expressed in terms of bases) is the matrix product of
a column vector with a row vector, and that this representation
is unique up to a choice of scale, so when we projectivize
(and thus introduce another choice of scale) we obtain
$$
X\simeq \BP A\times \BP B.
$$
$X$ is called the {\it Segre variety} and is often
written $X=Seg(\BP A\times \BP B)\subset \BP (A\ot B)$.

We calculate $\cC_{2,x}$ for the Segre. We first
must calculate $T_xX\subset T_x\BP V$. We identify
$T_x\BP V$ with $V\tmod \hat x$ and locate
$T_xX $ as a subspace of $V\tmod \hat x$.

Let $x=[a_0\ot b_0]\in Seg(\BP A\times \BP B)$. A curve $x(t)$ in $X$
with $x(0)=x$
is given by curves $a(t)\subset A$, $b(t)\subset B$, with
$a(0)=a_0,b(0)=b_0$ by taking $x(t)=[a(t)\ot b(t)$]. 
$$
\frac d{dt}|_{t=0}a_t\ot b_t= a_0'\ot b_0 + a_0\ot b_0'
$$
and thus
$$
T_xX = (A/a_0)\ot b_0 \op a_0\ot B/b_0 \tmod a_0\ot b_0
$$
Write $A'=(A/a_0)\ot b_0$, $B'=  a_0\ot (B/b_0)$ so
$$T_xX\simeq A'\op B'.
$$

We now take second derivatives modulo the tangent space to see which tangent directions
have lines osculating to order two (these will be the derivatives that
are zero modulo the tangent space).
\begin{align}
\frac {d^2}{(dt)^2}|_{t=0}a_t\ot b_t&=
a_0''\ot b_0 + a_0'\ot b_0' + a_0\ot b_0'' \tmod \hat x\\
&\equiv a_0'\ot b_0' \tmod \hat T_xX
\end{align}
Thus we get zero iff either $a_0'=0$ or $b_0'=0$, i.e.,
$$
\cC_{2,X,x}= \BP A'\sqcup \BP B'\subset \BP   (A'\op B')
$$
i.e., $\cC_{2,X,x}$ is the disjoint union of two linear spaces,
of dimensions $\tdim A-2, \tdim B-2$.
Note that $\tdim \cC_{2,x}$ is much
larger than expected.

For example, consider the case $Seg(\pp 2\times \pp 2)\subset \pp 8$.
Here $\cC_{2,x}$ is defined by four quadratic polynomials on
$\pp 3=\BP (T_xX)$, so one would have expected $\cC_{2,x}$ to be empty.
This rather extreme pathology led Griffiths and Harris to
conjecture:

\begin{conjecture}[Griffiths-Harris, 1979 \cite{GH}]
Let $Y^4\subset \pp 8$ be a variety not contained in  a hyperplane and let $y\in Y_{general}$.
If $\cC_{2,Y,y}=\pp 1\sqcup \pp 1\subset \pp 3=\BP (T_yY)$, then
$Y$ is isomorphic to $Seg(\pp 2\times \pp 2)$.
\end{conjecture}

(The original statement of the conjecture was in terms of the
projective second fundamental form defined below.)
Twenty years later, in \cite{Lrigid} I showed the conjecture was true, and moreover
in \cite{Lrigid,Lchss} I showed:

\begin{theorem}\label{chssrigid}
Let $X^n=G/P\subset \BP V$ be a rank two compact Hermitian symmetric
space (CHSS) in its minimal homogeneous embedding, other than a quadric
hypersurface. Let $Y^n\subset \BP V$ be a variety not contained
in a hyperplane and let $y\in Y_{general}$.
If $\cC_{2,Y,y}\simeq \cC_{2,X,x}$ then $Y$ is projectively isomorphic to $X$.
\end{theorem}

An analogous result is true in the $C^{\infty}$ category, namely

\begin{theorem}
Let $X^n=G/P\subset \BP V$ be a rank two compact Hermitian symmetric
space (CHSS) in its minimal homogeneous embedding, other than a quadric
hypersurface. Let $Y^n\subset \BP W$ be a smooth submanifold not contained
in a hyperplane.
If $\cC_{2,Y,y}\simeq \cC_{2,X,x}$  for all $y\in Y$, then $Y$ is projectively isomorphic to
an open subset of $X$.
\end{theorem}

The situation of the quadric hypersurface is explained below
(Fubini's theorem) - to characterize
it, one must have $\cC_{3,Y,y}=\cC_{3,Q,x}$.

The rank two CHSS are $Seg(\BP A\times \BP B)$, the
Grassmanians of two-planes $G(2,V)$, the quadric hypersurfaces, the
complexified Cayley plane
$\BO\BP^2= E_6/P_6$, and the
spinor variety $D_5/P_5$ (essentially the isotropic $5$-planes
through the origin in $\BC^{10}$ equipped with a quadratic form - the set of such planes is disconnected and the spinor variety is
one (of the two) isormorphic components. The minimal homogenous
embedding is also in a smaller linear space than the Plucker
embedding of the Grassmannian.)  The only rank one CHSS is projective space $\pp n$. 
The rank two CHSS and projective space are examples of 
{\it rational homogeneous varieties}.

\subsection{Homogeneous varieties}
Let $G\subset GL(V)$ be a reductive group acting irreducibly
on a vector space $V$. Then there exists a unique
closed orbit $X=G/P\subset \BP V$, which is called
a {\it rational homogeneous variety}.
(Equivalently, $X$ may be characterized as the orbit of a highest weight line,
or as the minimal orbit.)

Note that if $X=G/P\subset \BP V$ is homogenous, $T_xX$
inherits additional structure beyond that of a vector space.
Namely, consider $x=[Id]$ as the class of the identity
element for the projection $G\ra G/P$. Then $P$ acts on $T_xX$
and, as a $P$-module $T_xX\simeq \fg/\fp$. For example, in the
case of the Segre, $T_xX$ was the direct sum of two vector spaces.

A homogeneous variety $X=G/P$ is a {\it compact Hermitian symmetric
space}, or CHSS for short, if $P$ acts irreducibly on $T_xX$.
The {\it rank} of a CHSS is the  number of its last nonzero fundamental form
in its minimal homogeneous embedding. (This definition
agress with the standard one.) Fundamental forms are defined in \S\ref{ffsect}.
 
\begin{exercise}  The Grassmannian of $k$-planes through
the origin in $V$, which we denote $G(k,V)$, is homogenous
for $GL(V)$ (we have already seen the special case
$G(1,V)=\BP V$).

Determine the group $P_k\subset GL(V)$ that
stabilizes a point.
Show that $T_EG(k,V)\simeq E^*\ot V/E$ in two different
ways - by
an argument as in the Segre case above and by determining the structure of $\fg/\fp$.
\end{exercise}

While all homogeneous varieties have many special properties, the
rank at most two CHSS (other than the quadric hypersurface) are distinguished by the following property:

\begin{proposition}\cite{LM0} Theorem \ref{chssrigid} is sharp in the sense
that no other homogeneous variety is completely determined by
its asymptotic directions at a general point other than a linearly embedded
projective space.
\end{proposition}

Nevertheless, there are significant generalizations of theorem \ref{chssrigid}
due to Hwang-Yamaguchi and Robles discussed below in \S\ref{rigidsect}. To state these
results we will need definitions of the fundamental forms and Fubini cubic
forms, which are given in the next section.
 
\smallskip

However, with   an additional hypothesis - namely that the unknown
variety   has the correct codimension, we obtain the following result (which 
appears here for the first time):

\begin{theorem}\label{turbohy}
Let $X^n\subset \BC\pp\na$ be a complex submanifold not
contained in a hyperplane.  Let $x\in X$ be
a general point. Let  $Z^n\subset\pp{n+a}$ be an irreducible
  compact
Hermitian symmetric space in its minimal homogeneous embedding, other
than a quadric hypersurface.
 If $\cC_{2,X,x}=\cC_{2,Z,z}$   
 then $\overline X =Z$.
\end{theorem}

\begin{remark} 
The Segre variety has $\cC_{2, x}=\cC_{x}=\BP A'\sqcup \BP B'
\subset \BP (A'\op B')=\BP T_xX$.
 To see this, note
that a matrix has rank one iff all its $2\times 2$ minors
are zero, and these minors provide defining equations
for the Segre. In general, if a variety is defined by equations
of degree at most $d$, then any line having contact to order
$d$ at any point must be contained in the variety.  
In fact, by
an unpublished result of Kostant,  all homogeneously embedded rational homogeneous varieties $G/P$
are cut out by quadratic equations so $\cC_{2,G/P,x}=\cC_{G/P,x}$.
\end{remark}

\section{Moving frames and differential invariants}\label{movingframesect}
For more details regarding this section, see chapter 3 of \cite{EDSBK}. 

Once and for all fix index ranges $1\leq \a,\b,\g\leq n$, $n+1\leq
\mu,\nu\leq\na$, $0\leq A,B,C\leq n+a=N$.

\subsection{The Maurer-Cartan form of $GL(V)$}
Let $\tdim V=N+1$,
denote an element $f\in GL(V)$ by $f=(e_0\hd e_N)$ where we may think
of the $e_A$ as column vectors providing a basis of $V$. (Once a reference
basis of $V$ is fixed, $GL(V)$ is isomorphic to the space of all bases of $V$.) 
Each $e_A$ is a $V$-valued function on $GL(V)$,
$e_A: GL(V)\ra V$. For any differentiable map between manifolds
we can compute the induced differential
$$
de_A|_f: T_fGL(V)\ra T_{e_A}V
$$
but now since $V$ is a vector space, we may identify $T_{e_A}V\simeq V$ and
consider
$$
de_A: T_fGL(V)\ra V
$$
i.e., $de_A$ is a $V$-valued one-form on $GL(V)$.
As such, we may express it as
$$
de_A=e_0\ooo 0A + e_1\ooo 1A+\cdots + e_N\ooo NA
$$
where $\ooo AB\in \Omega^1(GL(V))$ are ordinary one-forms (This is
because $e_0\hd e_N$ is a basis of $V$ so any $V$-valued one
form is a linear combination of these with scalar valued one
forms as coefficients.) Collect the forms $\ooo AB$ into a matrix
$\Omega= (\ooo AB)$. Write $df=(de_0\hd de_N)$, so $df=f\Omega$ or
$$
\Omega= f\inv df
$$
$\Omega$ is called the {\it Maurer-Cartan form} for $GL(V)$.
Note that $\ooo AB$ measures the infinitesimal motion of
$e_B$ towards $e_A$.

\noindent{\bf Amazing fact}: we can compute the exterior derivative of $\Omega$ algebraically!
We have $d\Omega= d(f\inv)\ww df$ so we need to calculate
$d(f\inv)$. Here is where an extremely useful fact 
 comes in:
 
 \medskip

{\it The derivative of a constant function is zero.}

\medskip

We calculate $0=d(Id)=d(f\inv f)=d(f\inv)f+f\inv df$, and thus
$d\Omega=-f\inv df f\inv \ww df$ but we can move the scalar valued-matrix
$f\inv$ across the wedge product to conclude
$$
d\Omega = -\Omega\ww \Omega
$$ 
which is called the {\it Maurer-Cartan equation}.
The notation is such that $(\Omega \ww  \Omega)^A_B=\ooo AC\ww \ooo CB$.

\subsection{Moving frames for $X\subset \BP V$}
Now let $X^n\subset\pp{n+a}=\BP V$ be a submanifold. We are ready to take derivatives. Were we working in coordinates, to take
derivatives at $x\in X$, we might want to choose coordinates such that
$x$ is the origin. We will make the analoguous adaptation using moving frames,
but the advantage of moving frames is that all points will be as if they were the origin of a coordinate
system. To do this, let $\pi : \cf^0_X := GL(V)|_X\ra X$ be the restriction
of $\pi : GL(V)\ra \BP V$.
 
\smallskip

Similarly, we might want to choose local coordinates $(\xx 1\hd \xx{n+a})$
about $x=(0\hd 0)$ such that $T_xX$ is spanned by $\frac\partial{\partial x^1},\hd\frac\partial{\partial
x^n}$.  Again,   using moving frames the effect will be as if we had
chosen such coordinates about each point simultaneously.
To do this,  let $\pi :\cf^1\ra X$ denote the sub-bundle of
$\cf^0_X$ preserving the flag
$$
\hat x\subset\hat T_xX\subset V.
$$
Recall $\hat x\subset V$ denotes the line corresponding to $x$ 
and $\hat T_xX$ denotes the affine tangent space $T_v\hat X\subset V$, where
$[v]=x$.  Let $(\ee 0\hd \ee\na)$ be a
basis of $V$ with dual basis $(\ue 0\hd \ue\na)$ adapted such that
$\ee 0\in \hat x$ and $\{\ee 0,\ee \a\}$ span $\hat T_xX$. 
Write $T=T_xX$ and $N=N_xX=T_x\BP V/T_xX$.
 
\begin{remark}(Aside for the experts)
I am slightly abusing notation in this section by
identifying $\hat T_xX/\hat x$ with $T_xX:= (\hat T_xX/\hat x)\ot \hat x^*$
and similarly for $N_xX$.
\end{remark}

The fiber of $\pi : \cf^1\ra X$ over a point is isomorphic to the
group
$$
G_1=\left\{
\ g = \begin{pmatrix} g^0_0 & g^0_{\beta} & g^0_{\nu} \\
0 &  g^{\alpha}_{\beta} & g^{\alpha}_{\nu} \\
0 &  0 & g^{\mu}_{\nu} \end{pmatrix} \Big|\; g\in GL(V)  \right\}.
$$
 
While $\cf^1$ is not in general a Lie group, since $\cf^1\subset
GL(V)$, we may pull back the Maurer-Cartan from on $GL(V)$ to $\cf^1$.
Write the pullback of the Maurer-Cartan form to $\cf^1$ as
$$
\omega=\begin{pmatrix}\oo 0 & \ooo 0\beta & \ooo 0\nu\\
\oo\alpha & \ooo\alpha\beta & \ooo\alpha\nu\\
\ooo\mu 0 & \ooo\mu\beta & \ooo\mu\nu\end{pmatrix} .
$$
The definition of $\cf^1$ implies that $\oo\mu =0$
because 
$$de_0=\oo 0e_0+\cdots +\oo ne_n+\oo{n+1}\ee{n+1}+\cdots +\oo\na e_{n+1}
$$
but we have required that $e_0$ only move towards $e_1\hd e_n$ to
first order. Similarly, because $\tdim X=n$, the adaptation implies that
the forms $\oo\a$ are all linearly independent.

At this point you should know what to do - seeing something equal to
zero, we differentiate it. Thanks to the Maurer-Cartan equation, we may
calculate the derivative algebraically. We obtain
$$
0=d(\oo\mu)=-\ooo\mu\a\ww\oo\a \ \forall \mu
$$
Since the one-forms $\oo\a$ are all linearly independent,
  it is clear that the $\ooo\mu\a$ must be linear combinations
of the $\oo\b$, and in fact the {\it Cartan lemma} (see e.g.,  \cite{EDSBK}, p 314) implies that
the dependence is symmetric. More precisely (exercise!) there exist
functions
$$
\qq\mu\a\b : \cf^1\ra \BC
$$
with $\ooo\mu\a=\qq\mu\a\b\oo\b$ and moreover $\qq\mu\a\b=\qq\mu\b\a$.
One way to understand the equation $\ooo\mu\a=\qq\mu\a\b\oo\b$ is that the infinitesimal motion of the embedded tangent
space (the infinitesimal motion of the $e_{\a}$'s in the direction of
the $e_{\mu}$'s)    is determined by the motion of $e_0$ towards
the $e_{\a}$'s and the coeffiencients $\qq\mu\a\b$ encode
this dependence.

Now $\pi: \cf^1\ra X$ was defined   geometrically (i.e., without making
any arbitrary choices) so any function on $\cf^1$ invariant under the
action of $G_1$ descends to be a well defined function on $X$, and
will be a {\it differential invariant}. Our functions $\qq\mu\a\b$ are not
invariant under the action of $G_1$, but we can form a tensor
from them that is invariant, which will lead to a {\it vector-bundle
valued differential invariant} for $X$ (the same phenomenon happens
in the Euclidean geometry of submanifolds).

  Consider 
$$\tilde {II}_f=F_{2,f}:=\oo\a\ooo\mu\a\ot (\ee\mu \tmod \hat T_xX)=
\qq\mu\a\b\oo\a\oo\b\ot (\ee\mu \tmod \hat T_xX)
$$
$\tilde {II} \in \Gamma (\cf^1,\pi^*(S^2T^*X\ot NX))$ is constant on the fiber
and induces a tensor
$II\in \Gamma (X, S^2T^*X\ot NX)$. called the {\it projective second fundamental form}.

Thinking of $II_x: N^*_xX\ra S^2T^*_xX$, we may now properly define the asymptotic directions by
$$
\cC_{2,x}:= \BP (Zeros(II_x(N^*_xX))\subset \BP T_xX
$$

\subsection{Higher order differential invariants: the Fubini   forms}
We   continue differentiating constant functions:
$$
0=d(\ooo\mu\a -\qq\mu\a\b\oo\b)
$$
yields functions $\rr\mu\a\b\g:\cf^1\ra \BC$, symmetric
in their lower indices,  that induce a tensor
$F_3\in \Gamma (\cf^1,\pi^*(S^3T^*X\ot NX))$ called the
{\it Fubini cubic form}. Unlike the second fundamental form,
it does {\it not} descend   to be a tensor over $X$ because
it varies in the fiber. We discuss this variation   in 
\S\ref{hysect}. Such tensors provide {\it relative differential
invariants} and by sucessive differentiations, one
obtains a series of invariants 
$F_k\;\in\; \Gamma\;
(\cf^1,$ $\pi^*(S^kT^*\ot N))$. For example,
$$
\begin{aligned}
F_3&= \rr\mu\alpha\beta\gamma\oo\a\oo\b\oo
\gamma\ot\ee\mu \label{f3coeff}\\
F_4 &=  \rr\mu\alpha\beta{\gamma\delta}  
\oo\a \oo\b \oo\g \oo\delta\ot \ee\mu 
\end{aligned}
$$
where the functions $\rr\mu\alpha\beta\gamma,
\rr\mu\alpha\beta{\gamma\delta} $ are given by 
\begin{align}
\rr\mu\alpha\beta\gamma\oo\gamma 
&=-d\qq\mu\alpha\beta - \qq\mu\alpha\beta\ooo 00 
-\qq\nu\alpha\beta\ooo\mu\nu
+\qq\mu\alpha\delta\ooo\delta\beta +
 \qq\mu\beta\delta\ooo\delta\alpha \label{f3eqn}\\
\rr\mu\alpha\beta{\gamma\delta}\oo\delta  
&=  -d\rr\mu\alpha\beta\gamma -2\rr\mu\alpha\beta\gamma\ooo 00
-\rr\nu\alpha\beta\gamma\ooo\mu\nu \label{f4eqn}\\
&\quad+\mathfrak  S_{\alpha\beta\gamma}(
\rr\mu\alpha\beta\epsilon\ooo\epsilon\gamma 
+ \qq\mu\alpha\beta\ooo 0\gamma
-\qq\mu\alpha\epsilon\qq\nu\beta\gamma\ooo\epsilon\nu )\nonumber
\end{align}

We define $\cC_{k,x}:=Zeros(F_{2,f}\hd F_{k,f})\subset \BP T_xX$,
which is independent of our choice of $f\in \pi\inv(x)$.

If one chooses local affine coordinates $(\xx 1\hd \xx{n+a})$
such that $x=(0\hd 0)$ and $T_xX=\langle \frac\partial{\partial\xx\a}\rangle$, and writes $X$ as a graph
$$
\xx\mu = \qq\mu\a\b\xx\a\xx\b - \rr\mu\a\b\g\xx\a\xx\b\xx\g
+\rrr\mu\a\b\g\delta \xx\a\xx\b\xx\g\xx\delta + \cdots
$$
then there exists a  local section of $\cf^1$
such that 
$$
\begin{aligned}
F_2|_x&=\qq\mu\a\b d\xx\a d\xx\b\ot\frp{\xx\mu}\\
F_3|_x&=\rr\mu\a\b\g d\xx\a d\xx\b d\xx\g\ot\frp{\xx\mu}\\
F_4|_x&=\rrr\mu\a\b\g\delta d\xx\a d\xx\b d\xx\g 
d\xx\delta\ot\frp{\xx\mu}
\end{aligned}
$$
and similarly for higher orders.

\smallskip

\eqref{f3eqn} is a system of $a\binom{n+1}2$ equations
with one-forms as coefficients for the $a\binom{n+2}3$ coefficients
of $F_3$ and is overdetermined if we assume $d\qq\mu\a\b=0$, as we do
in the rigidity problems. One can calculate
directly  in the Segre
$Seg(\pp m\times \pp r)$, $m,r>1$ case that the only possible solutions are normalizable to zero
by a fiber motion as described in \eqref{hy3}. The situation is the same for $F_4,F_5$
in this case.

In general, once $F_k\hd F_{2k-1}$ are normalized to zero at a general point,
it is automatic that all higher $F_j$ are zero, see \cite{Lci}.
Thus one has the entire Taylor series and has completely identified the variety.
This was the method of proof used in \cite{Lrigid}, as the rank two
CHSS have all $F_k$ normalizable to zero when $k>2$.

Another perspective, for those familiar with $G$-structures, is that one
obtains rigidity by reducing   $\cf^1$ to a smaller bundle which is isomorphic to $G$, where the
homogeneous model is $G/P$.

Yet another perspective, for those familiar with exterior differential systems,
is that after three prolongations, the EDS defined by $I=\{\oo\mu,\ooo\mu\a-\qq\mu\a\b\oo\b\}$
on $GL(V)$ becomes involutive, in fact Frobenius.

\subsection{The higher fundamental forms}\label{ffsect}
A component of  $F_3$ does descend to a well defined tensor
on $X$. Namely, considering
$F_3: N^*\ra S^3T^*$, if we restrict $F_3|_{\tker F_2}$,
we obtain a tensor $\FF_3=III\in S^3T^*\ot N_3$ where $N_3=T_x\BP V/\{ T_xX + II(S^2T_xX)\}$.
One   continues
in this manner to get a series of tensors $\FF_k$ called the {\it fundamental forms}.

Geometrically, $II$ measures how $X$ is leaving its embedded tangent space at $x$ to
first order, $III$ measures how $X$ is leaving its second osculating space at
$x$ to first order while $F_3\tmod III$ measures how $X$ is moving away from
its embedded tangent space to second order.

\subsection{More rigidity theorems}\label{rigidsect}
Now that we have defined fundamental forms, we may state:

\begin{theorem} [Hwang-Yamaguchi]\cite{HY}
Let $X^n\subset \BC\pp\na$ be a complex submanifold.  Let $x\in X$ be
a general point. Let  $Z$ be an irreducible
rank $r$  compact
Hermitian symmetric space in its natural embedding, other
than a quadric hypersurface.
 If there exists linear maps $f: T_xX\ra T_zZ$, $g_k: N_{k,x}X\ra
N_{k,z}Z$ such that the induced maps
$S^kT^*_xX\ot N_xX\ra S^kT^*_zZ\ot N_zZ$ take $\BF_{k,X,x}$ to
$\BF_{k,Z,z}$   for $2\leq k\leq r$,
 then $\overline X =Z$.
\end{theorem}

In \cite{LM1} we calculated the differential invariants of the {\it adjoint
varieties}, the closed orbits in the projectivization of the adjoint representation
of a simple Lie algebra. (These are
  the homogeneous complex contact manifolds in their natural homogeneous
embedding.)  The adjoint varieties have $III=0$, but, in all cases but $v_2(\pp{2n-1})=C_n/P_1\subset \BP(\fc_n)=
\BP(S^2\BC^{2n})$ which we exclude from discussion in the remainder of this paragraph,
the invariants $F_3,F_4$ are   not normalizable to zero, even though
$\cC_{3,x}=\cC_{4,x}=\cC_{2,x}=\cC_x$. In a normalized frame $\cC_x$ is contained in a hyperplane $H$
and $F_4$ is the equation of the {\it tangential variety} of $\cC_x$ in $H$, where
the tangential variety $\t(X)\subset \BP V$ of an algebraic manifold $X\subset \BP V$ is the union of the points on the embedded tangent
lines ($\pp 1$'s) to the manifold. In this case the tangential variety is a hypersurface
in $H$, except for $\fg=\fa_n=\fsl_{n+1}$ which is discussed below. Moreover,
$F_3$ consists of the defining equations for the singular locus of $\t (\cC_x)$.
In \cite{LM1} we speculated that the varieties $X_{ad}$, with the exception of
$v_2(\pp{2n-1})$ (which is rigid to order three, see \cite{Lrigid}) would be rigid to order four, but not three, due to the nonvanishing
of $F_4$. Thus the following result came as a suprise to us:

\begin{theorem}[Robles, \cite{robles}] Let $X^{2(m-2)}\subset \BC\pp{m^2-2}$ be a complex submanifold.  Let $x\in X$ be
a general point. Let  $Z\subset\BP \fsl_m$ be the  adjoint variety.
 If   there exist  linear maps $f: T_xX\ra T_zZ$, $g: N_{x}X\ra
N_{z}Z$ such that the induced maps
$S^kT^*_xX\ot N_xX\ra S^kT^*_zZ\ot N_zZ$ take $ F_{k,X,x}$ to
$ F_{k,Z,z}$   for $k=2,3$,
 then $\overline X =Z$. 
\end{theorem}

Again, the corresponding result holds in the $C^{\infty}$ category.

The adjoint variety of $\fsl_m=\fsl(W)$ has the geometric interpretation of the variety of
flags of lines inside hyperplanes inside $W$, or equivalently as the traceless, rank one
matrices. It has $\cC_{2,z}$ the union of two disjoint linear spaces in a hyperplane
in $\BP T_zZ$. The quartic $F_4$ is the square of  a quadratic
equation (whose zero set contains the two linear spaces), and the cubics
in $F_3$ are the derivatives of this quartic, see \cite{LM1}, \S 6.

\subsection{The prolongation property and proof of  theorem \ref{turbohy}}
The  precise restrictions $II$ places on the $F_k$ in general is not  known at this time.
However,  there is
a strong restriction $II$ places on the higher fundamental forms
that dates back to Cartan.
We recall a definition from exterior differential systems:

\medskip

Let $U,W$ be vector spaces. Given a linear subspace
$A\subset S^kU^*\ot W$, define the {\it $j$-th prolongation} of
$A$ to be $A\up j:=(A\ot S^jU^*)\cap (S^{k+j}U^*\ot W)$. Thinking
of $A$ as a collection of $W$-valued homogeneous polynomials on $U$, the
$j$-th prolongation of $A$ is the set of all
homogeneous $W$-valued polynomials of degree
$k+j$ on $U$ with the property that all their $j$-th order partial derivatives
lie in $A$.

\begin{proposition}[Cartan \cite{C1919} p 377] Let $X^n\subset \pp{n+a}$, and let
$x\in X$ be a general point. Then
$\FF_{k,x}(N_k^*)\subseteq \FF_{2,x}(N_2^*)\up{k-2}$. (Here
$W$ is taken to be the trivial vector space $\BC$ and $U=T_xX$.)
\end{proposition}

\begin{proposition}\cite{LM0} Let $X=G/P\subset \BP V$ be
a CHSS in  its minimal homogeneous embedding. Then 
$\FF_{k,x}(N_k^*)= \FF_{2,x}(N_2^*)\up{k-2}$.
Moreover, the only nonzero components of the $F_k$ are
the fundamental forms.

The only homogeneous varieties having the property that
the only nonzero components of the $F_k$ are
the fundamental forms are the CHSS.
\end{proposition}

\begin{proof}[proof of \ref{turbohy}] The strict prolongation property
for CHSS in their minimal homogenous embedding implies that any variety
with the same second fundamental form at a general point as a CHSS in its minimal homogeneous
embedding can have codimension at most that of the corresponding CHSS, and
equality holds iff all the other fundamental forms are the prolongations
of the second.
\end{proof}

\section{Bertini type theorems and applications}

The results discussed so far dealt with homogeneous varieties. We now
broaden our study to various pathologies of the $\cC_{k,x}$.

Let $T$ be a vector space.  The classical {\it Bertini theorem} implies that
for a linear subspace $A\subset S^2T^*$, if $q\in A$ is such
that $\trank(q)\geq \trank(q')$ for all $q'\in A$, then
$u\in q_{sing}:=\{ v\in T\mid q(v,w)=0\ \forall w\in T\}$ implies $u\in
{\rm Zeros}\,(A):=\{ v\in T\mid Q(v,v)=0\ \forall Q\in A\}$.

\begin{theorem}[Mobile Bertini] \label{bertthm}
Let $X^n\subset \BP V$ be a complex manifold and let $x\in X$ be a
general point.
Let $q\in II(N^*_xX)$ be a generic quadric.  Then
  $q_{sing}$
is tangent to a linear space on   $\overline X$.
\end{theorem}

For generalizations and variations, see \cite{Lchss}.

The result holds in the $C^{\infty}$ category if one replaces a general point
by all points and that the linear space is contained in $X$ as long as $X$
continues (e.g., it is   contained in $X$ if $X$ is complete).

\begin{proof} 
Assume $v=  e_1\in q_{sing}$ and $q=q^{n+1}_{\a\b}\oo\a\oo\b$. Our hypotheses imply
$\qq{n+1} 1\b=0$ for all $\b$. Formula \eqref{f3eqn} reduces to
$$
\rr{n+1} 11\b\oo\b= -\qq{\mu} 11\ooo{n+1}\mu .
$$
If $q$ is generic we are   working on a reduction of $\cf^1$
where    the
$\ooo\mu\nu$ are independent of each other and independent of the semi-basic forms
(although the $\ooo\a\b$ will no longer be independent of the 
$\oo\a,\ooo\mu\nu$); thus the
coefficients on both sides of the equality are zero, proving both the
classical Bertini theorem and 
$v\in r_{sing}$ where $r$ is a generic cubic in $F_3(N^*)$. Then
using the formula for $F_4$ one obtains $v\in Zeros(F_3)$ and
$v\in s_{sing}$ where $s$ is a generic element of $F_4(N^*)$.
One then concludes by induction.
\end{proof}

\begin{remark} The mobile Bertini theorem
essentially  dates back to B. Segre \cite{bsegre}, and was
rediscovered in various forms in \cite{GH,Ein}. Its primarary use is in the
study of varieties  $X^n\subset \BP V$ with defective dual varieties
$X^*\subset \BP V^*$, where the dual variety of a smooth variety is the  
set of tangent hyperplanes to $X$, which is usually
a hypersurface. The point is that a generic quadric
in $II(N^*_{x,X})$ (with $x\in X_{general}$) is singular of rank $r$ 
iff $\tcodim X^*=n-r+1$. 
\end{remark}
 
 \begin{example}\label{p1pnex} Taking $X=Seg(\BP A\times \BP B)$ and keeping the notations
 of above,
 let $Y$ have the same second fundamental form of $X$ at
 a general point $y\in Y$ so we inherit an identification
 $T_yY\simeq A'\op B'$. If $\tdim B=b>a=\tdim A$, then  mobile Bertini implies that $\BP B'$ is actually
 tangent to a linear space on $Y$, because the maximum rank
 of a quadric is $a-1$. So at any point $[b]\in  \BP B'$
there is even a generic quadric  
 singular at $[b]$. (Of course the directions
 of $\BP A'$ are also tangent to lines on $Y$ because   the Segre
 is rigid.) 
 
\smallskip
 
While in \cite{Lrigid} I did not calculate the rigidity of
$Seg(\pp 1\times \pp n)\subset\BP (\BC^2\ot\BC^{n+1})$,   the rigidity
follows from the same calculations, however one must take   additional
derivatives to get the appropriate vanishing of the Fubini forms. However
there is also an elementary proof of rigidity in this case using the
mobile Bertini theorem   \cite{Lchss}. Given a variety $Y^{n+1}\subset \pp N$ such that
at a general point $y\in Y$, $\cC_{2,y}$ contains a $\pp{n-1}$, by
  \ref{bertthm} the resulting $n$-plane field
on $TY$ is integrable and thus $Y$ is
ruled by $\pp n$'s. Such a variety arises necessarily from a
curve in the Grassmannian $G(n+1,N+1)$ (as the union of the points on
the $\pp n$'s in the curve). But in order to also have the $\pp 0$ factor
in $\cC_{2,y}$, such a curve must be a line and thus $Y$ must be
the Segre.

\end{example}

The mobile Bertini theorem describes consequences of $\cC_{2,x}$ being pathological.
Here are some results when $\cC_{k,x}$ is pathological for $k>2$.

\begin{theorem}[Darboux] Let $X^2\subset\pp{2+a}$ be an analytic submanifold
  and let $x\in X_{general}$. If there exists a line $l$
having contact to order three with $X$ at $x$, then $l\subset \overline X$. In other words,
for surfaces in projective space, 
$$
\cC_{3,x}=\cC_x \ \forall x\in X_{general}.
$$
\end{theorem}

The $C^{\infty}$ analogue holds replacing general points by all points and
lines by line segments contained in $X$.

There are several  generalizations of this result in \cite{Lline,Lnfac}. Here is one of them:

\begin{theorem}\cite{Lnfac} \label{linethm}
 Let $X^n\subset\pp{n+1}$ be an analytic submanifold
  and let $x\in X_{general}$. If
 $\Sigma\subseteq \cC_{k,x}$ is an irreducible component
with $\tdim \Sigma >n-k$, then
$\Sigma\subset \cC_x$.
\end{theorem}

The $C^{\infty}$ analog holds with the by now obvious modifications.

The proof is similar to that of the mobile Bertini theorem.

\begin{exercise}One of my favorite problems to put on an undergraduate differential geometry
  exam is: Prove that   a surface in Euclidean three
space that has more than two lines passing through each point  is a plane
(i.e., has an infinite number of lines passing through). In \cite{MP}, Mezzetti
and Portelli showed that a $3$-fold having more than six lines
passing through a general point
must have an infinite number. Show that an $n$-fold having more than $n!$ lines passing
through a general point must have an infinite number passing through each point.
(See \cite{Lnfac} if you need help.)
\end{exercise}

\medskip

The rigidity of the quadric hypersurface is a classical result:

\begin{theorem}[Fubini] Let $X^n\subset\pp{n+1}$ be an analytic submanifold
or algebraic variety and let $x\in X_{general}$. Say
 $\cC_{3,x}=\cC_{2,x}$. Let $r=\trank\cC_{2,x}$. Then
 \begin{itemize}
 \item If $r>1$, then $X$ is a quadric hypersurface of rank $r$.
 \item If $r=1$, then $X$ has a one-dimensional Gauss image. In
 particular, it is ruled by $\pp{n-1}$'s.
 \end{itemize}
\end{theorem}

In all situations, the dimension of the Gauss image of $X$ is $r$,
see \cite{EDSBK}, \S 3.4.

One way to prove Fubini's theorem (assuming $r>1$) is to first use   mobile
Bertini to see that $X$ contains large linear spaces, then to note that
degree is invariant under linear section, so one can reduce to the case of
a surface. But then it is elementary to show that the only
analytic surface that
is doubly ruled by lines is the quadric surface.

Another way to prove Fubini's theorem (assuming $r>1$) is to use moving
frames to reduce the frame bundle to $O(n+2)$.

\medskip

{\it 
What can we say in higher codimension?
Consider codimension two. What are the varieties $X^n\subset \pp{n+2}$ such that
for general $x\in X$ we have 
$\cC_{3,x}=\cC_{2,x}$? 
}

\medskip

Note that there are two principal difficulties in codimension
two. First, in codimension one, having $\cC_{3,x}=\cC_{2,x}$ implies
that $F_3$ is normalizable to zero - this is no longer true in codimension
greater than one. Second, in codimension one, there is just one 
quadratic form in $II$, so its only invariant is its rank. In larger
codimension there are moduli spaces, although for pencils at least
there are normal forms, convienently given in \cite{HP}.

\smallskip

\noindent For examples, we have inherited from Fubini's theorem:

\smallskip

0. $\pp n\subset \pp{n+2}$ as a linear subspace

0'. $Q^n\subset\pp{n+1}\subset\pp{n+2}$ a quadric

0''. A variety with a one-dimensional Gauss image.

\smallskip

\noindent To these it is easy to see the following are also possible:

\smallskip

1. The (local) product of a curve with a variety with a one dimensional 
Gauss image.

2. The intersection of two quadric hypersurfaces.

3. A (local) product of a curve with a quadric hypersurface.

\smallskip

\noindent There is one more example we have already seen several times
in these lectures:

\smallskip

4. The Segre $Seg(\pp 1\times \pp 2)\subset \pp 5$ or a
cone over it.

\smallskip

\begin{theorem}[Codimension two Fubini]\cite{LR} Let $X^n\subset\pp{n+2}$ be an
analytic submanifold   and let $x\in X_{general}$.
If $\cC_{3,x}=\cC_{2,x}$ then $X$ is (an open subset of) one of
$0,0',0'',1-4$ above.
\end{theorem}

Here if one were to work over $\BR$, the corresponding result would be more
complicated as there are more normal forms for pencils of quadrics.

\section{Applications to algebraic geometry}

One nice aspect of algebraic geometry is that spaces parametrizing
algebraic varieties tend to also be algebraic varieties (or at 
least {\it stacks}, which is the algebraic geometer's version of
an orbifold). For example, let $Z^n\subset \pp{n+1}$ be a hypersurface.
The study of the images of holomorphic maps $f:\pp 1\ra Z$, called
{\it rational curves on $Z$} is of interest to algebraic geometers
and physicists. One can break this into a series of problems
based on the degree of $f(\pp 1)$ (that is, the number of
points in the intersection $f(\pp 1)\cap H$ where $H=\pp{n}$ is
a general hyperplane). When the degree is one, these are just
the lines on $Z$, and already here there are many open questions.
Let $\BF(Z)\subset \BG(\pp 1,\pp{n+1})=G(2,\BC^{n+2})$ denote
the variety of lines (i.e. linear $\pp 1$'s) on $Z$. 

Note that if the lines are distributed evenly on 
$Z$ and $z\in Z_{general}$, then $\tdim \BF(Z)=\tdim \cC_{z,Z}+ n-1$.  
 We always have 
$\tdim\BF(Z)\geq \tdim \cC_{z,Z}+ n-1$ (exercise!), so the microscopic geometry
bounds the macroscopic geometry.

\begin{example} Let $Z=\pp n$, then $\BF(Z)=G(2,n+1)$. In particular,
$\tdim\BF(Z) =2n-2$ and this is the largest possible
dimension, and $\pp n$ is the only variety with
$\tdim\BF(Z)= 2n-2$.
\end{example}

Which are the varieties $X^n\subset \BP V$
with $\tdim \BF(X)=2n-3$?  A classical theorem states that
in this case $X$ must be a quadric hypersurface.

Rogora, in \cite{Rogora}, classified all $X^n\subset\pp{n+a}$ with
$\tdim \BF(X)=2n-4$,  with the extra hypothesis $\tcodim X> 2$. The
only \lq\lq new\rq\rq\  example is the Grassmannian $G(2,5)$.  
(A classification in codimension one would be quite difficult.)
A corollary of the codimension two Fubini theorem is that   Rogora's theorem 
in codimension two is nearly proved - nearly 
and not completely because
one needs to add the extra hypothesis that $\cC_{2,x}$ has only
one component or that $\cC_{3,x}=\cC_{2,x}$.

\medskip

If one has extra information about $X$, one can say more. Say
$X$ is a hypersurface of degree $d$. Then it is easy to
show that $\tdim \BF(X)\geq 2n-1-d$.

\begin{conjecture}[Debarre, de Jong] If $X^n\subset\pp{n+1}$ is smooth and
$n>d=\tdeg (X)$, then $\tdim \BF(X)= 2n-1-d$.
\end{conjecture}

Without loss of generality it would be sufficient to prove the
conjecture when $n=d$ (slice by linear sections to reduce the dimension).
The conjecture is easy to show when $n=2$, it was proven by
Collino when $n=3$, by Debarre when $n=4$, and the proof of the $n=5$ case
was  was the 
PhD thesis of R. Beheshti \cite{beheshti}. Beheshti's thesis had three
ingredients, a general lemma (that $\BF(X)$ could not be uniruled by
rational curves), theorem \ref{linethm} above, and a case by case
argument. As a corollary of the codimension two Fubini theorem one
obtains a new proof of Beheshti's theorem eliminating the case
by case argument (but  there is a different case by case
argument buried in the proof of the codimension two Fubini theorem).
More importantly, the techniques should be useful in either proving
the theorem, or pointing to where one should look for potential
counter-examples for $n>5$.

\section{Moving frames proof of the Hwang-Yamaguchi theorem}\label{hysect}

The principle of calculation in \cite{Lchss}
was to use mobile Bertini theorems and the decomposition of the spaces $S^dT^*\ot N$ into irreducible $R$-modules, where 
$R\subset GL(T)\times GL(N)$ is the subgroup preserving
$II\in S^2T^*\ot N$. One can isolate where each $F_k$ can \lq\lq live\rq\rq\ as the
intersection of two vector spaces (one of which is $S^kT^*\ot N$, the other   is  $(\fg\upperp)_{k-3}$ defined below). Then, since $R$ acts on fibers, we can decompose
$S^kT^*\ot N$ and $(\fg\upperp)_{k-3}$ into $R$-modules and in order for a module to appear, it most be in both the vector
spaces. This combined with mobile Bertini theorems reduces the calculations to almost nothing.

  Hwang and Yamaguchi  use representation theory
in a more sophisticated way via
  a theory developed by Se-ashi \cite{seashi}. What follows is a proof of their result in the language of moving frames.

Let $Z=G/P\subset \BP W$ be a CHSS in its minimal homogeneous embedding.
For the moment we restrict to the case where $III_Z=0$ (i.e., $\trank Z=2$).
Let $X\subset \BP V$ be an analytic submanifold, let $x\in X$ be a general point
and assume $II_{X,x}\simeq II_{Z,z}$. We determine sufficient conditions
that imply $X$ is projectively isomorphic to $Z$.

  We have a filtration of $V$, $V_0=\hat x\subset V_1= \hat T_xX\subset V=V_2$.
Write $L=V_0,T=V_1/V_0,N=V/V_1$.
We have an induced grading of $\fgl(V)$ where 
\begin{align*}
\fgl(V)_0&=\fgl(L)\op \fgl(T)\op \fgl(N),\\
\fgl(V)_{-1}&=L^*\ot T \op T^*\ot N,\\
\fgl(V)_{1}&=L \ot T^* \op T \ot N^*,\\
\fgl(V)_{-2}&=L^*\ot N,\\
\fgl(V)_{2}&=L\ot N^*.
\end{align*}

In what follows we can no longer ignore the twist in defining $II$,
that is, we have  $II\in S^2T^*\ot N\ot L$.

Let $\fg_{-1}\subset \fgl(V)_{-1}$ denote the image of
\begin{align*}
T&\ra L^*\ot T + T^*\ot N\\
e_{\a}&\mapsto
e^0\ot e_{\a} + \qq\mu\a\b e^{\b}\ot e_{\mu}
\end{align*}

Let $\fg_0\subset \fgl (V)_0$ denote the subalgebra annhilating $II$. More precisely
$$
u=\begin{pmatrix}\xxx 00  & & \\
 & \xxx\a\b & \\
& & \xxx\mu\nu
\end{pmatrix}\in \fgl(V)_0
$$ 
is in $\fg_0$ iff
\begin{equation}\label{g0action}
u.II:= (-\xxx\mu\nu\qq\nu\a\b +\xxx\g\a\qq\mu\g\b + \xxx\g\b\qq\mu\a\g
-\xxx 00\qq\mu\a\b)e^{\a}\circ e^{\b}\ot e_{\mu}\ot e_0 =0.
\end{equation}

Let $\fg_1\subset \fgl(V)_1$ denote the maximal subspace such that
$[\fg_1,\fg_{-1}]=\fg_0$. Note that $\fg=\fg_{-1}\op \fg_0\op \fg_1$ coincides
with the $\BZ$-graded semi-simple Lie algebra giving rise to $Z=G/P$
and that the inclusion $\fg\subset \fgl(V)$ coincides with
the embedding $\fg\ra \fgl(W)$.  
The grading
on $\fg$ induced from the grading of $\fgl(V)$ agrees with the grading induced by $P$.
In particular,   $\fg_{\pm 2}=0$.
 
We let $\fg\upperp= \fgl(V)/\fg$   and note that
$\fg\upperp$ is naturally a $\fg$-module.  Alternatively,
one can work with $\fsl (V)$ instead of $\fgl(V)$ and define $\fg\upperp$ as 
the Killing-orthogonal complement to $\fg_0$ in   $\fsl(V)$ 
(Working with $\fsl(V)$-frames   does not effect projective geometry 
since the actual group of projective transformations
is $PGL(V)$). 

\smallskip

Recall that $F_3$ arises by applying the Cartan lemma to the equations
\begin{equation}\label{hyf3a}
0=-\ooo\mu\a\ww\ooo\a\b-\ooo\mu\nu\ww\ooo\nu\b +\qq\mu\a\b(\ooo\a 0\ww\oo 0 +
\ooo\a\g\ww\oo\g) \ \forall \mu,\b
\end{equation}
i.e., the tensor
\begin{equation}\label{hyf3}
(-\qq\nu\b\g \ooo\mu\nu +\qq\mu\a\g \ooo\a\b
 +
\qq\mu\a\b\ooo\a\g-\qq\mu\g\b \oo 0 )\ww\oo\g\ot  e^{\b}\ot e_{\mu}
\end{equation}
must vanish.
All forms appearing in the term in parenthesis in \eqref{hyf3}
are $\fgl(V)_0$-valued. Comparing with \eqref{g0action}, we see
that the $\fg_0$-valued part will
be zero and   $ (\fg\upperp)_0$ bijects to the image in
parenthesis.

Thus 
we may think of   obtaining the coefficients of $F_3$ at $x$ in two stages,
first we write the $(\fg\upperp)_0$ component of the Maurer-Cartan form
of $GL(V)$ as an arbitrary linear combination
of semi-basic forms, i.e., we choose a map
$T\ra(\fg\upperp)_0$. Once we have chosen such a map,
substituting the image into    \eqref{hyf3} yields a
$ (\La 2 T^*\ot \fgl(V)_{-1})$-valued tensor.
But by the definition of $\fg$, it is actually 
$ (\La 2 T^*\ot (\fg\upperp)_{-1})$-valued.
Then we require moreover that that this tensor is zero. 
In other words, pointwise we have a map
$$
\partial^{1,1}: T^*\ot (\fg\upperp)_0 \ra \La 2T^*\ot (\fg\upperp)_{-1}
$$
and the (at this stage) admissible coefficients of $F_3$ are determined by a choice of map $T\ra(\fg\upperp)_0$ which is in the kernel
of $\partial^{1,1}$.

\smallskip

Now the variation of $F_3$ as one moves in the fiber is given by
a map $T\ra (\fg\upperp)_0$ induced from the action of
$(\fg\upperp)_1$ on $ (\fg\upperp)_0$. We may express it as:
\begin{equation}\label{hy3}
\xxx 0\b e_0\ot e^{\b} + \xxx\a\nu e_{\a}\ot e^{\nu}
\mapsto 
\xxx 0\b\oo\b e_0\ot e^0+(\xxx 0\b\oo\a +\xxx\a\nu\ooo\nu\b)\ot e^{\b}\ot e_{\a}
+\xxx\a\nu\ooo\mu\a\ot e^{\nu}\ot e_{\mu}
\end{equation}
where 
$$
\begin{pmatrix}
0&\xxx 0\b &0\\ 0&0 &\xxx\g\nu \\ 0&0&0
\end{pmatrix}
$$
is a general element of $\fgl(V)_1$. 
The kernel of this map is    $\fg_1$ 
so it induces a linear map with source $(\fg\upperp)_1$. Similarly, the image automatically
takes values in $T^*\ot (\fg\upperp)_0$.
In summary, \eqref{hy3} may be expressed as a map
$$\partial^{2,0}: (\fg\upperp)_1\ra T^*\ot (\fg\upperp)_0.
$$

Let $C^{p,q}:=\La q T^*\ot (\fg\upperp)_{p-1}$. Considering $\fg\upperp$ as
a $T$-module (via the embedding $T \ra \fg$) we have the Lie algebra
cohomology groups
$$H^{p,1}(T, \fg\upperp):=
\frac{\tker \partial^{p,1}: C^{p,1}\ra C^{p-1,2}}
{\tim \partial^{p+1,0}: C^{p+1,0}\ra C^{p,1}}.
$$

 We summarize the above discussion:

\begin{proposition}
Let $X\subset \BP V$ be an analytic submanifold, let $x\in X_{general}$ and
suppose $II_{X,x}\simeq II_{Z,z}$ where $Z\subset \BP W$ is a
rank two CHSS in its
minimal homogeneous embedding.
The choices of $F_{3,X,x}$  imposed by \eqref{hyf3a}, modulo motions in the fiber, is isomorphic to the
Lie algebra cohomology group
$
H^{1,1}(T,\fg\upperp)
$.
\end{proposition}

Now if $F_3$ is normalizable and normalized to zero,
differentiating  again, we obtain
the equation

$$
(\qq\mu\a\b\ooo 0\g +\qq\mu\a\ep\qq\nu\b\g\ooo\ep\nu)\ww\oo\g
\ot e_{\mu}\ot e^{\a}\circ e^{\b}=0
$$
which determines the possible coefficients of $F_4$.

We conclude any choice of $F_4$ must be in the kernel of the map
$$
\partial^{2,1}: T^*\ot (\fg\upperp)_1\ra \La 2T^*\ot (\fg\upperp)_0
$$

The variation of  of $F_4$ in the fiber of $\cf^1$ is given by the image of the map
$$
\xxx 0\nu e_0\ot e^{\nu} \mapsto \xxx 0\nu\oo\ep e_{\ep}\ot e^{\nu}
$$
as $v\ww w$ ranges over the decomposable elements of $\La 2 T$. 
Without indicies, the variation of $F_4$ is the image of the map
$$
\partial^{3,0}: (\fg\upperp)_2\ra  T^*\ot (\fg\upperp)_1.
$$

We conclude that  if $F_3$ has been normalized to zero, then $F_4$ is normalizable to
zero if $H^{2,1}(T,\fg\upperp)=0$.

Finally, if $F_3,F_4$ are normalized to zero, then the coefficients of $F_5$ are
given by 
$$\tker\partial^{3,1}: T^*\ot (\fg\upperp)_2\ra \La 2T^*\ot (\fg\upperp)_1,
$$
and   there is nothing to quotient by in this case
because $\fgl(V)_3=0$.
In summary
\begin{proposition}\label{suffprop}
A sufficient condition for second order rigidity to hold for a rank two CHSS
in its minimal homogeneous embedding $Z=G/P\subset \BP W$
is that $H^{1,1}(T,\fg\upperp)=0$,$H^{2,1}(T,\fg\upperp)=0$,$H^{3,1}(T,\fg\upperp)=0$.
\end{proposition}

Assume for simplicity that $G$ is simple and $P=P_{\a_{i_0}}$
is the maximal parabolic subgroup obtained
by deleting all root spaces whose roots have a negative coefficient on
the simple root $\a_{i_0}$.
By Kostant's results \cite{kostant}, the $\fg_0$-module
$H^{*,1}(T,\Gamma)$ for {\it any} irreducible $\fg$-module $\Gamma$ of
highest weight $\l$
is the irreducible $\fg_0$-module with highest weight
$\s_{\a_{i_0}}(\l +\rho)-\rho$, where    $\s_{\a_{i_0}}$ is the   simple
reflection in the Weyl group
corresponding to $\a_{i_0}$ and $\rho$ is half the sum of the simple roots.
 
But now as long as $G/P_{\a_{i_0}}$ is not projective space or
a quadric hypersurface, 
Hwang and Yamaguchi, following \cite{SYY},  observe that
any non-trivial $\fg$-module $\Gamma$ yields
a $\fg_0$-module in $H^{*,1}$ that has a non-positive grading, so
one concludes that the above groups are {\it a priori} zero.
Thus one only need show that the trivial representation
is not a submodule of $\fg\upperp$.

\begin{remark} We note that the condition  in
proposition \ref{suffprop}  is {\it not} necessary for second
order rigidity. It holds in all rigid rank two cases but one,
$Seg(\pp 1\times \pp n)$, with $n>1$, which we saw, in \S\ref{p1pnex}, is
indeed rigid to order two. Note that in that case the na\"\i ve moving frames approach is
significantly more difficult as one must prolong several times before obtaining the
vanishing of the normalized $F_3$.
\end{remark}

To prove the general case of the Hwang-Yamaguchi theorem, 
say the last nonzero fundamental form is the $k$-th. Then one must show
$H^{1,1}\hd H^{k+1,1}$ are all zero. 
But again, Kostant's theory applies to show
all groups $H^{p,1}$ are zero for $p>0$.  Note that in this case, $H^{1,1}$
governs the vanishing   $F_{3,2}\hd F_{k,k-1}$, where $F_{k,l}$ denotes the
component of $F_k$ in $S^kT^*\ot N_l$.
In general, $H^{l.1}$ governs $F_{2+l,2}\hd F_{k+l,k-1}$.


\begin{thebibliography}{aa}

 
\bibitem{beheshti} R. Beheshti,   {\it Lines on projective hypersurfaces}, 
to appear in J. Reine Angew. Math.



\bibitem{bou} N. Bourbaki,
{\it Groupes et alg\`ebres de Lie}, 
Hermann, Paris, 1968,
MR0682756.

\bibitem{C1919} E. Cartan, {\it Sur les vari\'et\'es de courbure constante d'un
espace euclidien ou non euclidien},
{  Bull. Soc. Math France} 47 (1919) 125--160
and  48 (1920), 132--208; see also pp.~321--432 in
{  Oeuvres Compl\`etes} Part 3, Gauthier-Villars, 1955.


\bibitem{Ein} L. Ein,
{\it Varieties with small dual varieties, I},
{ Inventiones Math.} 86 (1986), 63--74.

\bibitem{fub} G. Fubini, 
{\it Studi relativi all'elemento lineare
proiettivo di una ipersuperficie}, 
Rend. Acad. Naz. dei Lincei, 1918, 99--106.  
 
\bibitem{GH} P.A. Griffiths \& J. Harris, 
{\it Algebraic Geometry and
Local Differential Geometry}, 
Ann. scient. Ec. Norm. Sup. {\bf 12} (1979) 355--432,
MR0559347.

\bibitem{Hwang} J.-M. Hwang,  
{\it  Geometry of minimal rational curves on Fano manoflds}, 
ICTP lecture notes, www.ictp.trieste.it./\~{}pub\_off/services.

\bibitem{HM} J.-M. Hwang \& N. Mok, 
{\it Uniruled projective manifolds with 
irreducible reductive $G$-structures}, 
J. Reine Angew. Math. {\bf 490} (1997) 55--64,
MR1468924.


\bibitem{HP}
W. V. D. Hodge \& D. Pedoe, 
{\it Methods of algebraic geometry. Vol. II.}   
Cambridge University Press, Cambridge (1994) p. 394+.

\bibitem{HY}
J.-M. Hwang \& K. Yamaguchi, 
{\it Characterization of Hermitian symmetric spaces by fundamental forms},
Duke Math. J.  120  (2003),  no. 3, 621--634. 

\bibitem{EDSBK} T. Ivey \& J.M. Landsberg,
{\it Cartan for beginners: differential geometry via moving frames and
exterior differential systems}, 
Graduate Studies in Mathematics, {\bf 61},
American Mathematical Society, Providence, RI, 2003,
MR2003610.


\bibitem{JM} G. Jensen  and E. Musso,
{\it Rigidity of hypersurfaces in complex projective space},
{  Ann. scient. Ec. Norm.} 27 (1994), 227--248.

\bibitem{kostant}   Kostant, Bertram 
{\it Lie algebra cohomology and the generalized Borel-Weil theorem}.  Ann. of Math. (2)  {\bf 74}  1961 329--387.
MR0142696

\bibitem{Lrr} J.M. Landsberg,
{\it On second fundamental forms of projective varieties},
Inventiones Math. {\bf 117} (1994) 303--315,
MR1273267.

\bibitem{Lci} \bysame,
{\it Differential-geometric characterizations of complete
intersections}, 
J. Differential Geom. {\bf 44} (1996) 32--73,
MR1420349.

\bibitem{Lrigid} \bysame, 
{\it On the infinitesimal rigidity of homogeneous varieties}, 
Compositio Math. {\bf 118} (1999) 189--201,
MR1713310.

\bibitem{Lsnu} \bysame,
{\it Algebraic geometry and projective differential geometry}, 
Seoul National University concentrated lecture series 1997, 
Seoul National University Press, 1999, 
MR1712467.

\bibitem{Lline} \bysame,
{\it Is a linear space contained in a submanifold?
- On the number of derivatives needed to tell},
{  J. reine angew. Math.} 508 (1999), 53--60.


\bibitem{Lnfac}\bysame,
{\it Lines on projective varieties}.  J. Reine Angew. Math.  {\bf 562}  (2003), 1--3. MR2011327

\bibitem{Lchss} \bysame,
{\it Griffiths-Harris rigidity of compact Hermitian symmetric spaces}, 
to appear in Journal of Differential 
Geometry (2006).

\bibitem{LM0} J.M. Landsberg \& L. Manivel, 
{\it On the projective geometry of rational homogeneous varieties},
Comment. Math. Helv. {\bf 78}({\bf 1}) (2003) 65--100,
MR1966752.

\bibitem{LM1} \bysame, 
{\it Classification of simple Lie algebras via projective geometry},  
Selecta Mathematica {\bf 8} (2002) 137--159,
MR1890196.


\bibitem{LMpop} \bysame,  {\sl Representation theory and projective 
geometry},  Algebraic Transformation Groups and Algebraic Varieties, 
Ed. V. L. Popov, Encyclopaedia of Mathematical Sciences {\bf 132}, Springer 2004,
 71-122. 


\bibitem{LMseries} \bysame,
{\it Series of Lie groups}, 
Michigan Math. J. {\bf 52}({\bf 2}) (2004) 453--479,
MR2069810.

\bibitem{LMdel} \bysame, {\sl  Triality, exceptional Lie algebras,
 and  Deligne dimension formulas}, Adv.  Math. {\bf 171} (2002), 59-85.


\bibitem{LMsex} \bysame,   {\it  The sextonions and $E\sb {7\frac 12}$}.  Adv. Math.  {\bf 201}  (2006),  no. 1, 143--179. MR2204753


\bibitem{LMuniv} \bysame,  {\it  A universal dimension formula for complex simple Lie algebras}.  Adv. Math.  {\bf 201}  (2006),  no. 2, 379--407. MR2211533


\bibitem{LMleg} \bysame, 
{\it Legendrian varieties}, 
  math.AG/0407279. To appear in Asian Math. J.

\bibitem{LR}  
Landsberg J.M., Robles C., 
{\it Fubini's theorem in codimension two}, 
preprint  math.AG/0509227

\bibitem{MP}   Mezzetti, E.; Portelli, D. {\it On threefolds covered by lines}  Abh. Math. Sem. Univ. Hamburg  70  (2000), 211--238. 
MR1809546

\bibitem{robles}
Robles, C., {\it Rigidity of the adjoint variety of $\fsl_n$}, preprint
math.DG/0608471.

\bibitem{Rogora}  
E. Rogora, {\it Varieties with many lines},  
Manuscripta Math.  82  (1994),  no. 2, 207--226.

\bibitem{SYY}
Sasaki,   Yamaguchi,  Yoshida, 
{\it On the rigidity of differential systems modelled on Hermitian symmetric spaces and disproofs of a conjecture concerning modular interpretations of configuration spaces.} in
CR-geometry and overdetermined systems (Osaka, 1994), 318--354,
Adv. Stud. Pure Math., 25,
Math. Soc. Japan, Tokyo, 1997. 

\bibitem{seashi}
Y. Se-Ashi,{\it 
On differential invariants of integrable finite type linear differential equations},
Hokkaido Math. J. {\bf 17} (1988), no. 2, 151--195.
MR0945853  


\bibitem{bsegre}
B. Segre, {\it Bertini forms and Hessian matrices},
{  J. London Math. Soc.} 26 (1951), 164--176.


\end{thebibliography}
\end{document}